\documentclass[a4paper]{article}

\usepackage{amsmath}
\usepackage[dvips]{epsfig}
\usepackage[latin1]{inputenc}

\newcommand{\for}{\operatorname{For}}
\newcommand{\crk}{\operatorname{crk}}
\newcommand{\rk}{\operatorname{rk}}
\newcommand{\zero}{\widehat{0}}
\newcommand{\one}{\widehat{1}}

\newtheorem{theo}{Theorem}[section]
\newtheorem{prop}[theo]{Proposition}
\newtheorem{lemma}[theo]{Lemma}
\newtheorem{coro}[theo]{Corollary}

\newenvironment{proof}{\begin{trivlist}\item{\bf{Proof.}}}
  {\hfill\rule{2mm}{2mm}\end{trivlist}}

\title{On intervals in some posets of forests}
\date{\today} \author{Frédéric Chapoton}

\begin{document}
\maketitle

\begin{abstract}
  We compute the characteristic polynomials of intervals in some
  posets of leaf-labeled forests of rooted binary trees.
\end{abstract}

\setcounter{section}{-1}

\section{Introduction}

The aim of this article is to study the poset $\for(I)$ attached to a
finite set $I$ which was introduced in \cite{bessel} in relation with
a Hopf operad of forests of binary trees. The underlying set of
$\for(I)$ is the set of leaf-labeled forests of rooted binary trees
with label set $I$. The main result is the following theorem.
\begin{theo}
  The characteristic polynomial of any interval in the poset $\for(I)$
  has only nonnegative integer roots.
\end{theo}

Furthermore, an explicit description of the roots is obtained for all
intervals. In particular, this gives simple product expressions for
all Möbius numbers.

\begin{figure}
  \begin{center}
    \leavevmode 
    \epsfig{file=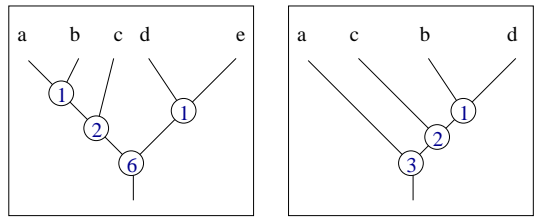,width=7cm} 
    \caption{Roots of characteristic polynomials.}
    \label{arbrexpo}
  \end{center}
\end{figure}

The simplest case is the interval between the minimal element $E$ of
the poset $\for(I)$ and a rooted binary leaf-labeled tree $T$ on $I$.
To each inner vertex of $T$, one associates the product of the number
of leaves of its two subtrees. These positive integers are the roots
of the characteristic polynomial of $[E,T]$. Figure \ref{arbrexpo}
displays two examples of this computation.

When the tree $T$ is a comb, the interval $[E,T]$ is isomorphic to the
partition lattice and the roots are $1,2,\dots,n$, where $n+1$ is the
cardinal of $I$, see the right example in Figure \ref{arbrexpo}. One
recovers the well-known factorization of the characteristic polynomial
of the partition lattice, by a method which differs from those
reviewed in \cite{sagan99}.

The other main result is an explanation of the coincidence of some
characteristic polynomials observed from the obtained description.
This is shown to be a consequence of some isomorphisms between the
intervals.

The strategy of proof is to decompose as much as possible the
intervals as products of simpler intervals. This gives a reduction to
the case of some special intervals, for which another kind of
decomposition can be done.

The first section is devoted to general results on these posets and to
the relation between combs and the partition lattice. The intervals
and their decompositions are studied in the second section. The third
section contains the proof that these posets are ranked by the number
of inner vertices. In the fourth section, invariants of the intervals
are computed, including the characteristic polynomials. The last
section contains the proof of some expected isomorphisms between the
intervals.

\section{Definition of posets}

\subsection{Notations}

A \textit{tree} is a leaf-labeled rooted binary tree and a
\textit{forest} is a set of such trees. Vertices are either inner
vertices (valence $3$) or leaves and roots (valence $1$). By
convention, edges are oriented towards the root. Leaves are
bijectively labeled by a finite set. Trees and forests are pictured
with their roots down and their leaves up, but are not to be
considered as planar. A leaf is an \textit{ancestor} of a vertex if
there is a path from the leaf to the root going through the vertex.

\smallskip

If $T_1$ and $T_2$ are trees on $I_1$ and $I_2$, let $T_1 \vee T_2$ be
the tree on $I_1 \sqcup I_2$ obtained by grafting the roots of $T_1$
and $T_2$ on a new inner vertex. If $F_1,F_2,\dots,F_k$ are forests on
$I_1,I_2,\dots,I_k$, let $F_1 \sqcup F_2 \sqcup \dots \sqcup F_k$ be
their disjoint union. If $F$ is the disjoint union of a forest on $J$
and a forest on $J'$, these restricted forests are denoted by $F[J]$
and $F[J']$. For a forest $F$, let $V(F)$ be the set of its inner
vertices. The number of trees in a forest $F$ on $I$ is the
difference between the cardinal of $I$ and the cardinal of $V(F)$.

\subsection{Posets of forests}

Let $F$ and $F'$ be forests on the label set $I$. Then set $F \leq F'$ if
there is a topological map from $F$ to $F'$ with the following
properties:
\begin{enumerate}
\item It is increasing with respect to orientation towards the root.
\item It maps inner vertices to inner vertices injectively.
\item It restricts to the identity of $I$ on leaves.
\item Its restriction to each tree of $F$ is injective.
\end{enumerate}

In fact, such a topological map from $F$ to $F'$ is determined up to
isotopy by the images of the inner vertices of $F$. One can recover
the map by joining the image of an inner vertex of $F$ in $F'$ with
the leaves of $F'$ which were its ancestor leaves in $F$.

Remark that there can be different $F$ lower than a given $F'$ with
the same image of $V(F)$ in $V(F')$.

\begin{lemma}
  \label{injectif}
  Let $F,F'$ be two distinct forests on $I$. If $F \leq F'$ then the
  cardinal of $V(F)$ is strictly less than the cardinal of $V(F')$.
\end{lemma}
\begin{proof}
  Assume that $F\leq F'$ and the cardinal of $V(F)$ is equal to that
  of $V(F')$. Then $F$ and $F'$ have the same number of trees. But
  each tree of $F$ is contained in a tree of $F'$ by connectivity.
  Each tree of $F'$ contains at least one tree of $F$ by injectivity
  on vertices. Therefore each tree of $F$ is contained in exactly one
  tree of $F'$. As these two trees have the same number of vertices,
  they must be equal. Hence $F=F'$.
\end{proof}

\begin{prop}
  The relation $\leq$ defines a partial order on the set $\for(I)$ of
  forests on $I$.
\end{prop}
\begin{proof}
  Reflexivity is given by the identity map. Transitivity is easy to
  check for each of the four required properties. Antisymmetry is
  clear by Lemma \ref{injectif}.
\end{proof}

\begin{figure}
  \begin{center}
    \leavevmode 
    \epsfig{file=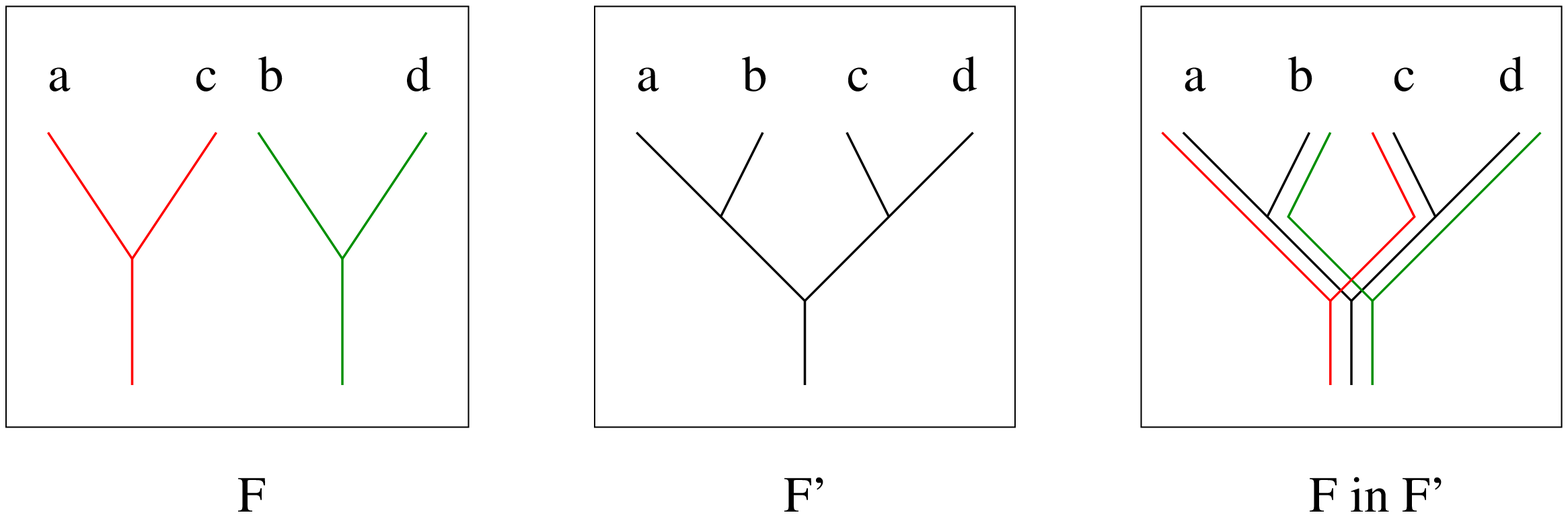,width=7cm} 
    \caption{A counterexample for the order relation.}
    \label{contrex}
  \end{center}
\end{figure}

\begin{figure}
  \begin{center}
    \leavevmode 
    \epsfig{file=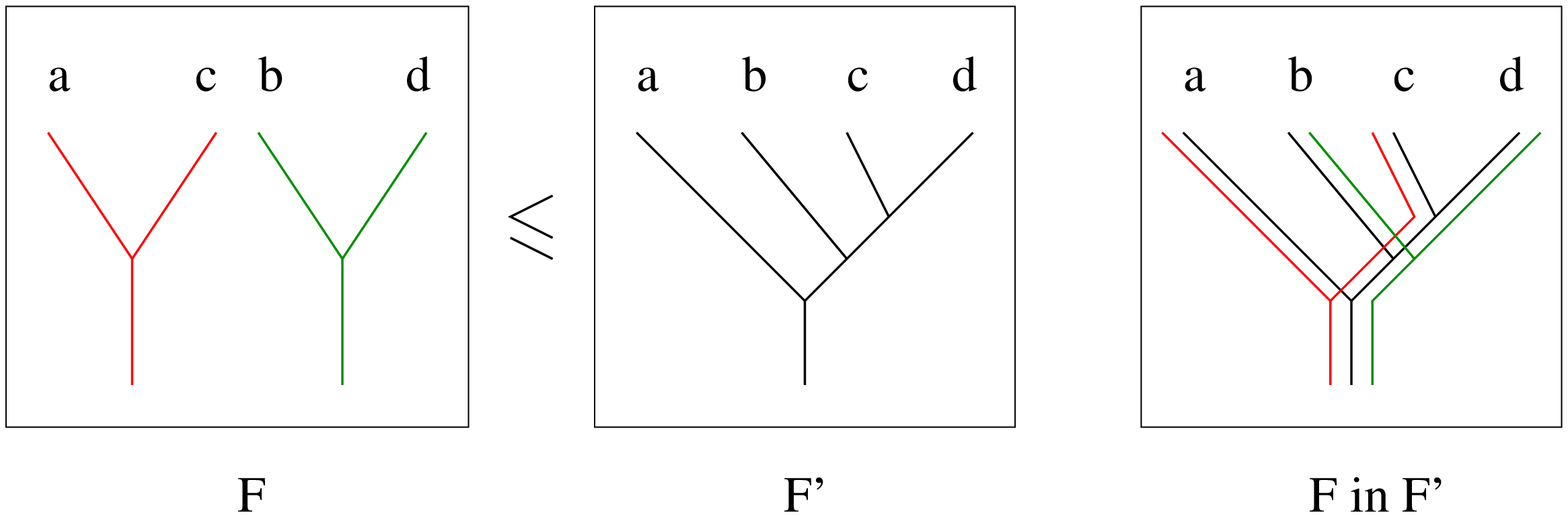,width=7cm} 
    \caption{An example for the order relation.}
    \label{exemple}
  \end{center}
\end{figure}

A counterexample, not injective on inner vertices, is given in Figure
\ref{contrex} and an example in Figure \ref{exemple}.

\begin{lemma}
  \label{2arbres}
  If $T_1$ and $T_2$ are trees on $I_1$ and $I_2$ then $T_1 \sqcup T_2
  \leq T_1 \vee T_2$.
\end{lemma}
\begin{proof}
  Obvious.
\end{proof}

Lemma \ref{2arbres} implies that, for each forest $F$ which is not a
tree, there exists a forest $F'$ with strictly less trees such that $F
\leq F'$. Lemma \ref{injectif} implies that trees are maximal
elements. Therefore the maximal elements of the poset $\for(I)$ are
exactly the trees. The forest without inner vertex is the unique
minimal element, denoted by $E$.

\begin{figure}
  \begin{center}
    \leavevmode 
    \epsfig{file=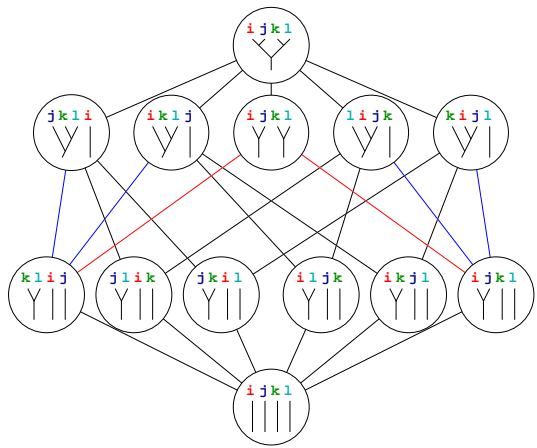,width=7cm} 
    \caption{An interval in the poset of forests on $\{i,j,k,\ell\}$.}
    \label{interval}
  \end{center}
\end{figure}

The intervals in the poset $\for(I)$ are not semimodular in general,
as can be seen on the interval depicted in Figure \ref{interval}.

\subsection{Relation to the partition lattice}

A \textit{comb} is a tree such that each inner vertex has at least one
of its two subtrees reduced to an edge.

\begin{prop}
  The interval between $E$ and a comb $C$ on the set $I$ is isomorphic
  to the partition lattice of the set $I$.
\end{prop}
\begin{proof}
  Remark first that a forest which is lower than a comb is necessarily
  composed of combs. The isomorphism $\phi$ is given by mapping a
  forest of combs to the partition of $I$ defined on the leaves by the
  combs. Let $J$ be a subset of $I$. Then there is exactly one comb
  $C_J$ with leaf set $J$ such that there exists an injective
  topological map from $C_J$ to $C$ which respects orientation and
  restricts to the identity of $J$ on leaves.
  
  This implies that each partition of $I$ can in only one way be
  realized as the leaf set of a forest of combs which is lower than
  $C$. Hence  $\phi$ is bijective. That the map $\phi$ is an
  isomorphism of posets follows easily from the description of the
  partial order, which is seen to coincide via $\phi$ with the
  refinement order on partitions.
\end{proof}

\section{Properties of intervals}

\subsection{Decomposition by connected components}

Let $F \leq F'$ be forests on $I$. Let $F'= T'_1,T'_2, \dots, T'_{k}$
seen as a set of trees $T'_j$ on $I_j$ with $I=I_1 \sqcup I_2 \sqcup
\dots \sqcup I_k$. Then $F$ can be uniquely decomposed as an union of
forests $F_j=F[I_j]$ on $I_j$ satisfying $F_j \leq T'_j$.

\begin{prop}
  \label{forestforest}
  The interval $[F,F']$ is isomorphic to the product of the intervals
  $[F_j,T'_j]$ in $\for(I_j)$ for $1 \leq j \leq k$.
\end{prop}
\begin{proof}
  Each element of this interval can in the same way be uniquely
  decomposed as an union of forests on $I_j$. The conditions defining
  the partial order then become equivalent to independent conditions
  on each part $I_j$.
\end{proof}

One can therefore restrict the attention to intervals between a forest
and a tree.

\subsection{Elements lower than a tree}

Let $T$ be a tree on the set $I$. Let us describe all elements $F$ of
$\for(I)$ which are lower than $T$. The binary tree $T$ defines a
partition $I=I_1 \sqcup I_2$ and two subtrees $T_1$ on $I_1$ and $T_2$
on $I_2$.

If $F_1 \leq T_1$ and $F_2 \leq T_2$ then clearly $F_1 \sqcup F_2 \leq
T$.

Let $F_1,F_2$ be forests with $F_1 \leq T_1$ and $F_2 \leq T_2$. Let
$J_1$ (resp. $J_2)$ be a chosen part of $I_1$ (resp. $I_2$)
corresponding to a chosen tree of $F_1$ (resp. $F_2$). Denote by
$G(F_1,J_1,F_2,J_2)$ the forest constructed from the disjoint union of
$F_1$ and $F_2$ by grafting a new inner vertex to the roots of the
chosen trees. This forest satisfies $G(F_1,J_1,F_2,J_2) \leq T$.

\begin{prop}
  \label{twokinds}
  A forest $F$ lower than $T$ is either the disjoint union $F_1 \sqcup
  F_2$ where $F_1\leq T_1$ and $F_2 \leq T_2$ or is equal to
  $G(F_1,J_1,F_2,J_2)$ where $F_1\leq T_1$, $F_2 \leq T_2$ and $J_1$,
  $J_2$ are parts of $I_1$,$I_2$ corresponding to trees of
  $F_1$,$F_2$.
\end{prop}

\begin{proof}
  Two forests $F_1$ and $F_2$ can be defined as follows. Consider the
  inner vertices of $F$ having only elements of $I_1$ as ancestors. By
  joining them in $T$ to their ancestor leaves, one gets $F_1$ on
  $I_1$ which satisfies $F_1 \leq T_1$. The same construction gives
  $F_2$ on $I_2$ with $F_2 \leq T_2$.
  
  Assume first that the image of $V(F)$ in $V(T)$ does not contain the
  lowest inner vertex of $T$. From the definition of the poset, $F$ is
  in fact lower than $T_1 \sqcup T_2$ and is the disjoint union $F_1
  \sqcup F_2$.
  
  Assume now on the contrary that the image of $V(F)$ in $V(T)$
  contains the lowest inner vertex of $T$. By injectivity on inner
  vertices, there exists a unique tree $T'$ of $F$ which has an inner
  vertex mapped to the lower inner vertex of $T$. By injectivity on
  trees, the tree $T'$ can be written $T'_1 \vee T'_2$ where $T'_1$
  (resp. $T'_2$) has leaf set $J_1 \subset I_1$ (resp. $J_2 \subset
  I_2$). The tree $T'_1$ (resp. $T'_2$) is a tree of $F_1$ (resp.
  $F_2$) and $F$ is indeed equal to $G(F_1,J_1,F_2,J_2)$.
\end{proof}




\subsection{Intervals under a tree}

\label{undertree}

Let $T$ be a tree and $F$ a forest on the set $I$ such that $F \leq T$
and the image of $V(F)$ in $V(T)$ contains the lowest inner vertex of
$T$. This implies that $F$ can be written $G(F_1,J_1,F_2,J_2)$ as
explained in the previous section.

\begin{prop}
  \label{markedtree}
  The interval $[F,T]$ is isomorphic to the product of the intervals
  $[F_1,T_1]$ in $\for(I_1)$ and $[F_2,T_2]$ in $\for(I_2)$.
\end{prop}
\begin{proof}
  Let $F'$ be an element of the interval $[F,T]$. Necessarily the
  image of $V(F')$ contains the lowest vertex of $T$. Therefore one
  can write $F'=G(F'_1,J'_1,F'_2,J'_2)$ with $F'_1 \leq T_1$ and $F'_2
  \leq T_2$. By definition of the partial order, the inequality $F
  \leq F'$ implies that $F_1\leq F'_1$, $F_2 \leq F'_2$ and that
  $J'_1$ (resp. $J'_2$) must contain $J_1$ (resp. $J_2$). It follows
  that $J'_1$ and $J'_2$ are uniquely determined for a given $F'_1$
  and $F'_2$. Therefore, any pair $(F'_1,F'_2)$ with $F_1\leq F'_1
  \leq T_1$ and $F_2 \leq F'_2 \leq T_2$ can be uniquely extended to
  an element of $[F,T]$.
  
  The elements of the interval $[F,T]$ are therefore in bijection with
  pairs $(F'_1,F'_2)$ in $[F_1,T_1]\times [F_2,T_2]$.
  
  The conditions defining the partial order do not depend on $J_1$ and
  $J_2$, and are mapped by the bijection to independent conditions on
  $I_1$ and $I_2$. Hence the bijection is an isomorphism of posets.
\end{proof}

\subsection{Special intervals}

\label{intree}

Let $F$ be a forest and $T$ be a tree on the set $I$ with $F\leq T$.
Assume that the image of $V(F)$ in $V(T)$ does not contain the lowest
inner vertex of $T$, that is to say $F$ is a disjoint union $F_1
\sqcup F_2$ on $I_1$ and $I_2$. The intervals of the form $[F,T]$ for
such $F$ and $T$ are called \textit{special intervals}.

\smallskip

\begin{prop}
  \label{threekinds}
  There are three kinds of sub-intervals in a special interval
  $[F,T]$:
  \begin{enumerate}
  \item $[F'_1 \sqcup F'_2,F''_1 \sqcup F''_2]$ with $F_1 \leq F'_1
    \leq F''_1 \leq T_1$ and $F_2 \leq F'_2 \leq F''_2 \leq T_2$.
    This interval is isomorphic to $[F'_1,F''_1] \times [F'_2,F''_2]$.
  \item $[G(F'_1,J'_1,F'_2,J'_2),G(F''_1,J''_1,F''_2,J''_2)]$ with
    $F_1 \leq F'_1 \leq F''_1 \leq T_1$ and $F_2 \leq F'_2 \leq F''_2
    \leq T_2$ where $J''_1$ and $J''_2$ are the unique parts of
    $F''_1$ and $F''_2$ containing $J'_1$ and $J'_2$. This interval is
    isomorphic to $[F'_1,F''_1] \times [F'_2,F''_2]$.
  \item $[F'_1 \sqcup F'_2,G(F''_1,J''_1,F''_2,J''_2)]$ with $F_1 \leq
    F'_1 \leq F''_1\leq T_1$, $F_2 \leq F'_2 \leq F''_2\leq T_2$, and
    $J''_1,J''_2$ are arbitrary parts of $F''_1$ and $F''_2$.
  \end{enumerate}
\end{prop}
\begin{proof}
  First, let us determine which elements $F'$ can be lower than $T$
  and greater than $F$. If $F'$ is a disjoint union $F'_1 \sqcup
  F'_2$, then it is necessary and sufficient that $F_1 \leq F'_1 \leq
  T_1$ and $F_2 \leq F'_2 \leq T_2$. If $F'=G(F'_1,J'_1,F'_2,J'_2)$,
  then necessary and sufficient conditions are also that $F_1 \leq
  F'_1 \leq T_1$ and $F_2 \leq F'_2 \leq T_2$.
  
  Let us discuss now the possible intervals according to the type of
  their bounds. First, it is not possible to have a relation
  $G(F'_1,J'_1,F'_2,J'_2) \leq (F''_1 \sqcup F''_2)$, because the
  lowest inner vertex is present in the first element and not in the
  second one, which would contradict injectivity.
 
  Let us study each of the three remaining cases.
  
  Case $[\sqcup,\sqcup]$: One can apply Prop. \ref{forestforest}. The
  interval $[F'_1 \sqcup F'_2,F''_1 \sqcup F''_2]$ is non-empty if and
  only if $F'_1 \leq F''_1$ and $F'_2 \leq F''_2$. If these conditions
  are fulfilled, this interval is isomorphic to the claimed product.
  
  Case $[G,G]$: the interval
  $[G(F'_1,J'_1,F'_2,J'_2),G(F''_1,J''_1,F''_2,J''_2)]$. Let
  $F''_1=f''_1 \sqcup T''_1$ where $T''_1=F''_1[J''_1]$ is a tree and
  similarly let $F''_2=f''_2 \sqcup T''_2$ where $T''_2=F''_2[J''_2]$
  is a tree. Then $G(F''_1,J''_1,F''_2,J''_2)$ is equal to $f''_1
  \sqcup f''_2 \sqcup (T''_1 \vee T''_2)$. One can then decompose the
  interval as a product by Prop. \ref{forestforest}. By applying Prop.
  \ref{twokinds} to the interval under $T''_1 \vee T''_2$, the product
  interval is non-empty if and only if one has $F'_1 \leq F''_1$ and
  $F'_2 \leq F''_2$ and the parts $J''_1$ and $J''_2$ contains
  respectively the parts $J'_1$ and $J'_2$. When these conditions are
  satisfied, Prop. \ref{markedtree} shows that this interval is
  isomorphic to the claimed product.
 
  Case $[\sqcup,G]$: the interval $[F'_1 \sqcup
  F'_2,G(F''_1,J''_1,F''_2,J''_2)]$. Let $F''_1=f''_1 \sqcup T''_1$
  where $T''_1=F''_1[J''_1]$ is a tree and similarly let $F''_2=f''_2
  \sqcup T''_2$ where $T''_2=F''_2[J''_2]$ is a tree. Then
  $G(F''_1,J''_1,F''_2,J''_2)$ is equal to $f''_1 \sqcup f''_2 \sqcup
  (T''_1 \vee T''_2)$. One can then decompose the interval as a
  product by Prop. \ref{forestforest}. By applying Prop.
  \ref{twokinds} to the interval under $T''_1 \vee T''_2$, the product
  interval is non-empty if and only if one has $F'_1 \leq F''_1$ and
  $F'_2 \leq F''_2$.


\end{proof}

\section{Rank property}

Say that a finite poset is \textit{ranked} if it has a unique minimal
element $\zero$ and all maximal chains have the same length. Note that
this definition differs slightly from the usual definition which
requires the uniqueness of the maximal element.

\begin{prop}
  The poset $\for(I)$ is ranked by the number of inner vertices.
\end{prop}
\begin{proof}
  The proof is by recursion on the cardinal of $I$. The proposition is
  true by inspection for small $I$.
  
  Fix a maximal interval $[E,T]$ where $T$ is a tree on $I$ and $E$ is
  the forest without inner vertices. Consider a maximal chain $E=F_0
  \leq \dots \leq F=F_{k-1} \leq F_k =T$ in $[E,T]$. It is clear from
  Lemma \ref{injectif} that the length $k$ is at most the number of
  inner vertices of $T$.
  
  Let us discuss according to the properties of $F$.
  
  Assume first that $F$ contains the lowest inner vertex of $T$. By
  maximality, there should be no element between $F$ and $T$, and one
  can conclude by recursion hypothesis and Prop. \ref{markedtree} that
  either $F_1=T_1$ and $F_2$ has just one vertex less than $T_2$ or
  the similar situation obtained by exchanging $1$ and $2$ holds.
  
  Assume on the contrary that $F$ does not contain the lowest inner
  vertex of $T$. By maximality, there should be no element between $F$
  and $T$, and one can conclude by recursion hypothesis and Prop.
  \ref{threekinds} that $F_1=T_1$ and $F_2=T_2$.
  
  Therefore, in both cases, the number of inner vertices of $F$ is the
  number of inner vertices of $T$ minus one.

  By recursion and Prop. \ref{forestforest}, the length of all maximal
  chains of $[E,F]$ is the number of inner vertices of $F$.
  
  This implies that the length of all maximal chains of $[E,T]$ is the
  number of inner vertices of $T$. All trees on $I$ have the same
  number of inner vertices. The proposition is proved.
\end{proof}

Note that the corank function in $\for(I)$ is given by the number of
trees minus one.

\section{Invariants of intervals}

For a standard reference on posets, see \cite{stanley1}.

\subsection{$M$-polynomials and $Z$-polynomials}

Let $P$ be a ranked poset with unique minimal element $\zero$ and
unique maximal element $\one$. Let $\crk$ be the corank function on
$P$, which is defined by $\crk(a)=\rk(\one)-\rk(a)$. The degree of the
poset is $\deg(P)=\crk(\zero)$ .

One defines the $M$-polynomial of the poset $P$, which is a generating
function for the Möbius function, as follows:
\begin{equation}
  M(P)=\sum_{a \leq b}\mu(a,b) x^{\crk(a)}y^{\crk(b)}.
\end{equation}

In the same way, one defines the $Z$-polynomial, which is a generating
function for the zeta function, as follows:
\begin{equation}
  Z(P)=\sum_{a \leq b}x^{\crk(a)}y^{\crk(b)}.
\end{equation}

The characteristic polynomial is defined to be
\begin{equation}
  \chi(P)=\sum_{b}\mu(\zero,b) y^{\crk(b)}.
\end{equation}

The cardinal polynomial is the generating function for the corank:
\begin{equation}
  \operatorname{Card}(P)=\sum_{a} x^{\crk(a)}.
\end{equation}

The Möbius number is $\mu(P)=\mu(\zero,\one)$.

It is clear that the $M$-polynomial allows to recover the
characteristic polynomial, which in turn contains the Möbius number as
leading coefficient. The $M$-polynomial also contains the information
of the cardinal polynomial.

The cardinal polynomial is also determined by the $Z$-polynomial.

The following proposition is classical.

\begin{prop}
  Let $P_1$ and $P_2$ be two such ranked posets and $P_1 \times P_2$
  their product. Then
  \begin{equation}
    Z_{P_1 \times P_2}=Z_{P_1}Z_{P_2}\text{  and  }
    M_{P_1\times  P_2}=M_{P_1}M_{P_2}.
  \end{equation}
\end{prop}

\begin{lemma}
  \label{value1}
  The value at $y=1$ of the $M$-polynomial is $1$.
\end{lemma}
\begin{proof}
  This is an immediate consequence of the definition of the Möbius
  function and the existence of $\one$.
\end{proof}

\subsection{$Z$-polynomials of special intervals}

Let $F$ be a forest and $T$ be a tree on the set $I$. Assume that $F
\leq T$ and the image of $V(F)$ in $V(T)$ does not contain the lowest
inner vertex of $T$. We keep the notations of section \ref{intree}.

\begin{theo}
  \label{vee-zeta}
  The $Z$-polynomial of the special interval $[F,T]$ (denoted by Z) is
  determined by the $Z$-polynomials of the intervals $[F_1,T_1]$ and
  $[F_2,T_2]$ (denoted by $Z_1$ and $Z_2$). One has
  \begin{equation}
    \label{Zvee}
    Z=xy Z_1 Z_2 + \partial_x (x Z_1)\partial_x (x Z_2)
    + x \partial_y (y Z_1)\partial_y (y Z_2).
  \end{equation}
\end{theo}

\begin{proof}
  The sum defining $Z$ is split in three parts, according to the three
  different kinds of subintervals in $[F,T]$ listed in Prop.
  \ref{threekinds}.
  
  The first part is given by
  \begin{equation*}
    \sum_{{F'_1 \leq F''_1} \atop
      {F'_2 \leq F''_2}}x^{1+\crk(F'_1)+\crk(F'_2)}
    y^{1+\crk(F''_1)+\crk(F''_2)},
  \end{equation*}
  which is $xy Z_1 Z_2$.

  The second part is given by
  \begin{multline*}
    \sum_{{F'_1 \leq F''_1} \atop
      {F'_2 \leq F''_2}}\sum_{J'_1,J'_2} x^{\crk(F'_1)+\crk(F'_2)}
    y^{\crk(F''_1)+\crk(F''_2)}\\=
    \sum_{{F'_1 \leq F''_1} \atop
      {F'_2 \leq F''_2}}(1+\crk(F'_1))(1+\crk(F'_2))x^{\crk(F'_1)+\crk(F'_2)}
    y^{\crk(F''_1)+\crk(F''_2)},
  \end{multline*}
  which is $\partial_x (x Z_1)\partial_x (x Z_2)$.

  The last part is given by
  \begin{multline*}
    \sum_{{F'_1 \leq F''_1} \atop
      {F'_2 \leq F''_2}}\sum_{J''_1,J''_2} x^{1+\crk(F'_1)+\crk(F'_2)}
    y^{\crk(F''_1)+\crk(F''_2)}\\=
    \sum_{{F'_1 \leq F''_1} \atop
      {F'_2 \leq F''_2}}(1+\crk(F''_1))(1+\crk(F''_2))
    x^{1+\crk(F'_1)+\crk(F'_2)}
    y^{\crk(F''_1)+\crk(F''_2)},
  \end{multline*}
  which is $x \partial_y (y Z_1)\partial_y (y Z_2)$.

  This concludes the proof of the theorem.
\end{proof}

\subsection{$M$-polynomials of special intervals}

Let $F$ be a forest and $T$ be a tree on the set $I$. Assume that $F
\leq T$ and the image of $V(F)$ in $V(T)$ does not contain the lowest
inner vertex of $T$. We keep the notations of section \ref{intree}.

\begin{theo}
  \label{vee-mobius}
  The $M$-polynomial of the special interval $[F,T]$ (denoted by $M$)
  depends only on the $M$-polynomials of the intervals $[F_1,T_1]$ and
  $[F_2,T_2]$ (denoted by $M_1$ and $M_2$). One has
  \begin{equation}
    \label{Mvee}
    M=xy M_1 M_2 + (1-x) \partial_x (x M_1)\partial_x (x M_2).
  \end{equation}
\end{theo}

\begin{proof}
  By recursion on the degree of $[F,T]$.
  
  Formula (\ref{Mvee}) is correct if $F=T_1 \sqcup T_2$ and $T=T_1 \vee
  T_2$, which is the only possible case of degree $1$.
 
  The sum which defines $M$ is split in three parts, according to the
  three different kinds of subintervals in $[F,T]$ listed in Prop.
  \ref{threekinds}.

  The first part is given by
  \begin{equation*}
    \sum_{{F'_1 \leq F''_1} \atop
      {F'_2 \leq F''_2}}\mu(F'_1,F''_1)\mu(F'_2,F''_2)x^{1+\crk(F'_1)+\crk(F'_2)}
    y^{1+\crk(F''_1)+\crk(F''_2)},
  \end{equation*}
  which is $xy M_1 M_2$.

  The second part is given by
  \begin{multline*}
    \sum_{{F'_1 \leq F''_1} \atop
      {F'_2 \leq F''_2}}\sum_{{J'_1}
      \atop {J'_2}} \mu(F'_1,F''_1)\mu(F'_2,F''_2)
    x^{\crk(F'_1)+\crk(F'_2)}
    y^{\crk(F''_1)+\crk(F''_2)}\\=
    \sum_{{F'_1 \leq F''_1} \atop
      {F'_2 \leq F''_2}}(1+\crk(F'_1))(1+\crk(F'_2))
    \mu(F'_1,F''_1)\mu(F'_2,F''_2)
    \\ x^{\crk(F'_1)+\crk(F'_2)}
    y^{\crk(F''_1)+\crk(F''_2)},
  \end{multline*}
  which is $\partial_x (x M_1)\partial_x (x M_2)$.
  
  The computation of the third part is more complicated. It is given
  by
  \begin{equation}
    \label{troisieme}
     \sum_{{F'_1 \leq F''_1} \atop
      {F'_2 \leq F''_2}}
   \bigg{(}   \sum_{{J''_1}\atop {J''_2}}
      \mu([F'_1 \sqcup F'_2,G(F''_1,J''_1,F''_2,J''_2)]) \bigg{)}
      x^{1+\crk(F'_1)+\crk(F'_2)}
    y^{\crk(F''_1)+\crk(F''_2)}.
  \end{equation}
  
  Let us simplify the inner summation. Fix a part $J''_1$ and a part
  $J''_2$ and assume that $\crk(G(F''_1,J''_1,F''_2,J''_2))\not=0$.
  This excludes only the case when $G(F''_1,J''_1,F''_2,J''_2)=T$.
  
  Then $G(F''_1,J''_1,F''_2,J''_2)$ is not a tree. Let $K''_1$ be the
  complement of $J''_1$ in $I_1$ and $K''_2$ be the complement of
  $J''_2$ in $I_2$.
  
  Denote the trees $F''_1 [J''_1]$ and $F''_2 [J''_2]$ by $T''_1$ and
  $T''_2$. Then $G$ can be uniquely decomposed as
  \begin{equation*}
    F''_1[K''_1] \sqcup F''_2[K''_2]  \sqcup (T''_1 \vee T''_2)
  \end{equation*}
  and $F'_1 \sqcup F'_2$ can also be decomposed as
  \begin{equation*}
    F'_1[K''_1] \sqcup F'_2[K''_2]  \sqcup (F'_1[J''_1] \sqcup F'_2[J''_2]).
  \end{equation*}
  
  Hence the interval $[F'_1 \sqcup F'_2,G(F''_1,J''_1,F''_2,J''_2)]$
  is isomorphic to the product
  \begin{equation*}
    [F'_1[K''_1],F''_1[K''_1]] \times [F'_2[K''_2],F''_2[K''_2]]
    \times [F'_1[J''_1] \sqcup F'_2[J''_2],T''_1 \vee T''_2].
  \end{equation*}
  
  Choose enumerations $J''_1(1),J''_1(2),\dots,J''_1(k_1)$ and \-
  $J''_2(1),J''_2(2),\dots,J''_2(k_2)$ of the parts of parts of
  $F''_1$ and $F''_2$. Remark that $k_1=\crk(F''_1)+1$ and
  $k_2=\crk(F''_2)+1$. The inner sum can be restated as
  \begin{multline*}
    \sum_{j=1}^{k_1} \sum_{\ell=1}^{k_2}
    \mu([F'_1[K''_1(j)],F''_1[K''_1(j)]])
    \mu([F'_2[K''_2(\ell)],F''_2[K''_2(\ell)]])\\
    \mu([F'_1[J''_1(j)] \sqcup F'_2[J''_2(\ell)],
    T''_1(j) \vee T''_2(\ell)]),
  \end{multline*}
  where $T''_1(j)$ and $T''_2(\ell)$ are obvious notations.
  
  As stated below in Corollary \ref{vee-carac}, it follows from the
  recursion hypothesis that for all special intervals
  $[F^\natural,T^\natural]$ of smaller degree, one has
  \begin{equation}
    \mu([F^\natural,T^\natural])=-(\deg([F^\natural_1,T^\natural_1])+1)
    (\deg([F^\natural_2,T^\natural_2])+1)
    \mu([F^\natural_1,T^\natural_1])\mu([F^\natural_2,T^\natural_2]).
  \end{equation}
  
  Hence using this consequence of the recursion hypothesis, the inner
  sum is
  \begin{multline*}
    -\sum_{j=1}^{k_1} \sum_{\ell=1}^{k_2}
    \mu([F'_1[K''_1(j)],F''_1[K''_1(j)]])
    \mu([F'_2[K''_2(\ell)],F''_2[K''_2(\ell)]]) \\ 
    \big((\deg([F'_1[J''_1(j)],T''_1(j)])+1)
    (\deg([F'_2[J''_2(\ell)],T''_2(\ell)])+1)\\
    \mu([F'_1[J''_1(j)],T''_1(j)])\mu([ F'_2[J''_2(\ell)],
    T''_2(\ell)])\big),
  \end{multline*}
  which in turn is equal to 
  \begin{equation}
    -(\crk(F'_1)+1)(\crk(F'_2)+1)\mu([F'_1,F''_1])\mu([F'_2,F''_2]).
  \end{equation}
  
  Hence the third part (\ref{troisieme}) is equal, up to a polynomial
  in $x$ corresponding to special intervals with maximal element $T$,
  to
  \begin{multline}
     -\sum_{{F'_1 \leq F''_1} \atop
      {F'_2 \leq F''_2}}
     (1+\crk(F'_1))(1+\crk(F'_2))
      \mu([F'_1,F''_1])\mu([F'_2,F''_2])\\
      x^{1+\crk(F'_1)+\crk(F'_2)}
    y^{\crk(F''_1)+\crk(F''_2)},
  \end{multline}
  which is $-x \partial_x (x M_1)\partial_x (x M_2)$.
  
  Therefore the full sum $M$ is equal to the expected formula, up to a
  polynomial in $x$. By Lemma \ref{value1} the value of $M$ at $y=1$
  is $1$, and the value of the right-hand-side of Formula (\ref{Mvee})
  at $y=1$ is also $1$. Hence Formula (\ref{Mvee}) stands exactly. The
  recursion step is done and the theorem is proved.
\end{proof}

\smallskip

\begin{coro}
  \label{vee-carac}
  The characteristic polynomial of the interval $[F,T]$ (denoted by
  $\chi$) depends only on the characteristic polynomials of the intervals
  $[F_1,T_1]$ and $[F_2,T_2]$ (denoted by $\chi_1$ and $\chi_2$). One has
  \begin{equation}
    \chi(y)= (y - (\deg_1+1)(\deg_2+1)) \chi_1(y) \chi_2(y),
  \end{equation}
  where $\deg_1=\deg([F_1,T_1])$ and $\deg_2=\deg([F_2,T_2])$. As a
  special case, one has 
  \begin{equation}
    \mu= -(\deg_1+1)(\deg_2+1)\mu_1 \mu_2
  \end{equation}
  for Möbius numbers.
\end{coro}



\subsection{Factorization of characteristic polynomials}

Let $F,F'$ be forests on the set $I$ with $F \leq F'$. Let $V$ be the
image of $V(F)$ in $V(F')$. Let us call \textit{marked vertices} the
elements of $V$. To each non-marked vertex $v \in V(F') \setminus V$,
there correspond two subtrees $T_1$ and $T_2$. Let $d_1$ (resp. $d_2$)
be the number of leaves of $T_1$ (resp. $T_2$) minus the number of
marked vertices of $T_1$ (resp. $T_2$). One associates to the
non-marked vertex $v$ its \textit{exponent} which is the integer $d_1
d_2$. The exponents of the pair $(F,F')$ are the exponents of the
non-marked vertices of $F'$.

Remark that the exponents of $(F,F')$ only depend on $F'$ and the set
$V$ of marked inner vertices, not on which $F$ is mapped to $F'$ using
$V$.

\begin{theo}
  The characteristic polynomial of the interval $[F,F']$ has only
  positive integer roots, which are the exponents of the pair
  $(F,F')$.
\end{theo}
\begin{proof}
  The proof is by recursion on the cardinal of $I$ and $\deg([F,F'])$.
  
  The statement is true for $\deg([F,F'])=0$ or $I$ a singleton, in
  which cases the set of exponents is empty.
  
  If $F'$ is not a tree, then using Prop. \ref{forestforest}, the
  statement is a consequence of the recursion hypothesis.
  
  If $F'$ is a tree and $V$ contains the bottom vertex of $F'$, then
  the statement follows from the recursion hypothesis by using Prop.
  \ref{markedtree}.
  
  If $F'$ is a tree and $V$ does not contain the bottom vertex of
  $F'$, then the statement follows from the recursion hypothesis by
  using Corollary \ref{vee-carac}.
  
  The theorem is proved.
\end{proof}

\begin{figure}
  \begin{center}
    \leavevmode 
    \epsfig{file=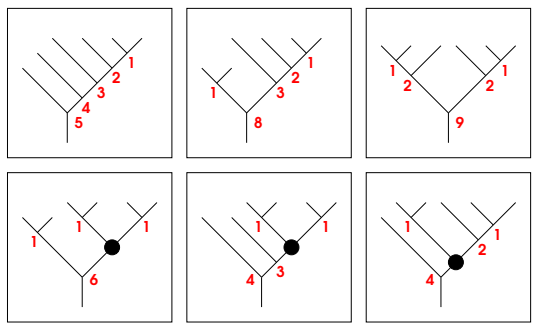,width=7cm} 
    \caption{Other examples of exponents.}
    \label{expofig}
  \end{center}
\end{figure}

Some more examples are given in Figure \ref{expofig}.

\section{Partitive posets}

\subsection{Definition and product}

A \textit{partitive poset} is a ranked poset $P$ with $\zero$ and
$\one$ together with a ranked-poset map $f_P$ from $P$ to a subposet
of the partition lattice of a finite set $I$ (possibly with shifted
rank function).

The only examples which will be used here are all intervals $[F,F']$
where $F$ and $F'$ are forests on $I$ with the map defined by the
partition of the label set according to trees.

One can define the product of two partitive posets $(P_1,I_1,f_1)$ and
$(P_2,I_2,f_2)$. As a ranked poset, it is the usual product $P_1
\times P_2$. The composition of $f_1 \times f_2$ with the inclusion
map of partition lattices (induced by disjoint union of partitions)
defines a poset map from $P_1 \times P_2$ to the partition lattice of
$I_1 \sqcup I_2$.

\begin{prop}
  \label{product}
  The interval $[F_1 \sqcup F_2,F'_1 \sqcup F'_2]$ is isomorphic
  as a partitive poset to the product of $[F_1,F'_1]$ and
  $[F_2,F'_2]$.
\end{prop}

\begin{proof}
  This is an easy reformulation of Prop. \ref{forestforest}.
\end{proof}

\subsection{Twisted product of partitive posets}

Consider two partitive posets $P_1,P_2$ such that their respective
$\one$ are mapped to a partition with only one part.

Choose for each of these posets a part of the image of their $\zero$
in the corresponding partition lattice. These chosen parts are denoted
by $K_1$ and $K_2$.

The twisted product of $P_1$ and $P_2$ is defined as follows. As a
poset, it is simply $P_1 \times P_2$. The map to a partition lattice
differs from the map for the usual product by gathering the parts
containing $K_1$ and $K_2$ to a single part.

It is easy to see that, up to isomorphism of partitive poset, this
construction does not depend on the choices made.

\smallskip

Let $F$ be a forest and $T$ be a tree on the set $I$ with $F \leq T$.
The binary tree $T$ defines a partition $I=I_1 \sqcup I_2$ and two
subtrees $T_1$ on $I_1$ and $T_2$ on $I_2$. Assume that the image of
$V(F)$ contains the lowest inner vertex of $T$. Recall the notations
of \ref{undertree}.

\begin{prop}
  \label{twisted}
  The interval $[F,T]$ is isomorphic as a partitive poset to the
  twisted product of the intervals $[F_1,T_1]$ and $[F_2,T_2]$.
\end{prop}

\begin{proof}
  This is an easy reformulation of Prop. \ref{markedtree}.
\end{proof}

\subsection{The $\vee$-product of partitive posets}

Consider two partitive posets $P_1,P_2$ such that their respective
$\one$ are mapped to a partition with only one part.


The $\vee$ product of $P_1$ and $P_2$ is defined as follows. The
underlying set is the disjoint union of $P_1 \times P_2$ with elements
denoted by $\{a_1 \sqcup a_2\}$ and of the set
$\{G(a_1,J_1,a_2,J_2)\}$ where $a_1$ and $a_2$ are elements of $P_1$
and $P_2$ respectively and $J_1$ and $J_2$ are parts of the image of
$a_1$ and $a_2$.



The order relation is given by the following relations: 
\begin{enumerate}
\item $(a_1 \sqcup a_2) \leq (a'_1 \sqcup a'_2)$ if $a_1 \leq a'_1$
  and $a_2 \leq a'_2$.
\item $G(a_1,J_1,a_2,J_2) \leq G(a'_1,J'_1,a'_2,J'_2)$ if $a_1 \leq
  a'_1$, $a_2 \leq a'_2$ and the part $J'_1$ (resp. $J'_2$) contains
  the part $J_1$ (resp. $J_2$).
\item $(a_1 \sqcup a_2) \leq G(a'_1,J'_1,a'_2,J'_2)$ if $a_1 \leq
  a'_1$ and $a_2 \leq a'_2$.
\end{enumerate}

The map to a partition lattice is defined as follows. An element $a_1
\sqcup a_2$ is mapped to the disjoint union of the partitions
associated to $a_1$ and $a_2$. An element $G(a_1,J_1,a_2,J_2)$ is
mapped to the partition obtained from the disjoint union of the
partitions associated to $a_1$ and $a_2$ by gathering the two parts
containing $J_1$ and $J_2$ to a single part.

The result is a partitive poset, called the $\vee$-product of $P_1$ and $P_2$.

\smallskip

Let $F$ be a forest and $T$ be a tree on the set $I$. Assume that $F
\leq T$ and the image of $V(F)$ in $V(T)$ does not contain the lowest
inner vertex of $T$. We keep the notations of section \ref{intree}.

\begin{prop}
  \label{veeprod}
  The interval $[F,T]$ is isomorphic as a partitive poset to the
  $\vee$-product of the intervals $[F_1,T_1]$ and $[F_2,T_2]$ .
\end{prop}

\begin{proof}
  This is essentially a reformulation of Prop. \ref{threekinds}.
\end{proof}

\subsection{Marked trees}

Let $F,F'$ be forests on the set $I$ with $F\leq F'$.

\begin{theo}
  Up to isomorphism of partitive posets, the interval $[F,F']$ depends
  only on the pair $(F',V)$ where $V$ is the subset of marked inner
  vertices of $F'$ associated to $F$.
\end{theo}
\begin{proof}
  By recursion on the degree of $[F,F']$ and the cardinal of $I$. This
  is clear if the degree is zero or the cardinal of $I$ is one.
  
  If $F'$ is not a tree, then the proposition follows from the
  recursion hypothesis and Prop. \ref{product}.
  
  If $F'$ is a tree and $V$ contains the lowest inner vertex of $F'$,
  the statement follows from the recursion hypothesis and Prop.
  \ref{twisted}.
  
  If $F'$ is a tree and $V$ does not contain the lowest inner vertex
  of $F'$, this follows from the recursion hypothesis and Prop.
  \ref{veeprod}.

\end{proof}

\nocite{*}

\bibliographystyle{plain}
\bibliography{treillis}
\end{document}